\documentclass{article}
\hoffset - 20 mm

\catcode`\@=11 \@addtoreset{equation}{section}

\catcode`\@=12

\newtheorem{Theorem}{Theorem}[section]
\newtheorem{Proposition}{Proposition}[section]
\newtheorem{Lemma}{Lemma}[section]
\newtheorem{Corollary}{Corollary}[section]

\newcommand{\bTheorem}[1]{
\begin{Theorem} \label{T#1} }
\newcommand{\eT}{\end{Theorem}}

\newcommand{\bProposition}[1]{
\begin{Proposition} \label{P#1}}
\newcommand{\eP}{\end{Proposition}}

\newcommand{\bLemma}[1]{
\begin{Lemma} \label{L#1} }
\newcommand{\eL}{\end{Lemma}}

\newcommand{\bCorollary}[1]{
\begin{Corollary} \label{C#1} }
\newcommand{\eC}{\end{Corollary}}

\newcommand{\bFormula}[1]{
\begin{equation} \label{#1}}
\newcommand{\eF}{\end{equation}}

\newcommand{\Ov}[1]{\overline{#1}}
\newcommand{\DC}{C^\infty_c}

\newcommand{\nep}{n_{\ep}}

\newcommand{\vue}{\vu_\ep}

\newcommand{\vu}{\vc{u}}
\newcommand{\vc}[1]{{\bf #1}}
\newcommand{\qed}{\bigskip \rightline {Q.E.D.} \bigskip}
\newcommand{\Div}{{\rm div}_x}
\newcommand{\Grad}{\nabla_x}

\newcommand{\tn}[1]{\mbox {\F #1}}
\newcommand{\dx}{{\rm d} {x}}
\newcommand{\dt}{{\rm d} t }

\newcommand{\intO}[1]{\int_{\Omega} #1 \ \dx}

\newcommand{\bProof}{{\bf Proof: }}

\newcommand{\ep}{\varepsilon}

\font\F=msbm10 scaled 1000

\date{}
\makeindex
\begin{document}

%%%%%%%%%%%%%%%%%%%%%%%%%%%%%%%%

%%%%%%%%%%%%%%%%%%%%%%%%%%%%%%%%

\title{On the vanishing electron-mass limit in plasma hydrodynamics in
unbounded media}
\author{Donatella Donatelli  \and Eduard Feireisl\thanks{The work of E.F. was supported by Grant 201/09/0917 of GA \v CR as a part of the general research
programme of the Academy of Sciences of the Czech Republic,
Institutional Research Plan AV0Z10190503.} \and Anton{\' \i}n Novotn\' y \thanks{The work of A.N. was partially supported by the general research
programme of the Academy of Sciences of the Czech Republic,
Institutional Research Plan AV0Z10190503.} }
\maketitle

\bigskip

\centerline{Dipartimento di Matematica Pura ed Applicata}
\centerline{Universit\` a degli Studi dell'Aquila, 67 100 L'Aquila, Italy}

\medskip

\centerline{Institute of Mathematics of the Academy of Sciences of the Czech Republic}
\centerline{\v Zitn\' a 25, 115 67 Praha 1, Czech Republic}

\medskip

\centerline{Universit\' e du Sud Toulon-Var, BP 20132, 839 57 La
Garde, France}

\medskip

\begin{abstract}

We consider the zero-electron-mass limit for the Navier-Stokes-Poisson system in unbounded spatial domains. Assuming smallness of the viscosity coefficient and ill-prepared initial data, we show that the asymptotic limit
is represented by the incompressible Navier-Stokes system,
with a Brinkman damping, in the case when viscosity is proportional to the electron-mass, and by the incompressible
Euler system provided the viscosity is dominated by the electron mass. The proof is based on the RAGE theorem and dispersive estimates for acoustic waves, and on the concept of suitable weak solutions for the compressible Navier-Stokes system.

\end{abstract}

%\tableofcontents

\section{Introduction}
\label{i}

Singular limits arise frequently in the process of \emph{model reduction} in fluid mechanics. In this paper we consider the limit of vanishing
ratio electron mass/ion mass in a hydrodynamic model for plasma confined to
an unbounded spatial domain $\Omega \subset R^3$.

\subsection{Equations}

For a given (constant) density $N_i$ of positively charged ions, the time evolution of the electron density $n_e = n_e(t,x)$ and the electron velocity $\vu = \vu(t,x)$ is governed by the system of equations
\bFormula{i1}
\partial_t n_e + \Div (n_e \vu)=0,
\eF
\bFormula{i2}
m_e \Big( \partial_t (n_e \vu) + \Div( n_e \vu \otimes \vu) \Big) +
\Grad p(n_e) = \Div \tn{S} (\Grad \vu) + n_e \Grad \Phi - m_e \frac{n_e \vu}{\tau},
\eF
\bFormula{i3}
\Delta \Phi= n_e - N_i,
\eF
where $m_e$ is the ratio of the electron/ions mass, $p$ is the electron
pressure, $\Phi$ is the electric potential, $\tau$ is the relaxation time, and $\tn{S}$ denotes the viscous stress tensor,
\bFormula{i4}
\tn{S}(\Grad \vu) = \mu \Big( \Grad \vu + \Grad^t \vu -\frac{2}{3} \Div \vu \tn{I} \Big),\ \mu > 0,
\eF
see Anile and Pennisi \cite{AniPen}, J\" ungel and Peng \cite{JunPen2},
\cite{JunPen1}.
Moreover, we suppose the electron velocity satisfies the slip boundary conditions
\bFormula{i5}
\vu \cdot \vc{n}|_{\partial \Omega} = 0, [\tn{S} (\Grad \vu) \vc{n}]
\times \vc{n}|_{\partial \Omega} = 0,
\eF
and the boundary is electrically insulated,
\bFormula{i6}
\Grad \Phi \cdot \vc{n}|_{\partial \Omega} = 0.
\eF
As the underlying spatial domain is unbounded, we also prescribe
prescribe the far field behavior:
\bFormula{i7}
\vu(x) \to 0 , \ \Phi(x) \to 0  \ \mbox{as}\ |x| \to \infty.
\eF

Our goal is to study the singular limit and identify the limit
problem for $m_e \to 0$ under the
condition
\begin{itemize}
\item
$\mu \approx m_e$, or
\item
$\mu / m_e \to 0$ as $m_e \to 0$.
\end{itemize}

\subsection{Ill-prepared initial data}

For $m_e = \ep^2$, $\mu_\ep = \mu/\ep^2$, problem (\ref{i1} - \ref{i7}) is reminiscent of the low Mach (incompressible) limit of the Navier-Stokes system that have been investigated in a number of recent studies, see
see the survey papers by
Danchin \cite{DANC3}, Gallagher \cite{Gallag}, Masmoudi \cite{MAS1},
and Schochet \cite{SCH2}, and the references cited therein. The zero-electron-mass limit for the \emph{inviscid} fluid was treated recetly
by Ali et al. \cite{AliCheJun}, see also Chen, Chen and Zhang  \cite{CheCheZha}. In the latter case, it is shown that the systems becomes neutral, meaning
$n_e \to N_i$, while the limit velocity field $\vc{v}$ satisfies a
damped Euler system
\bFormula{euler1}
\Div \vc{v} = 0,
\eF
\bFormula{euler2}
\partial_t \vc{v} + \vc{v} \cdot \Grad \vc{v} + \frac{1}{\tau} \vc{v} +
\Grad \Pi = 0,
\eF
supplemented with the impermeability condition
\bFormula{euler3}
\vc{v} \cdot \vc{n}|_{\partial \Omega} = 0.
\eF

In \cite{AliCheJun}, \cite{CheCheZha}, the authors consider the periodic boundary conditions and \emph{well-prepared} initial data
\[
n_e(0, \cdot) = N_i + \ep^2 N_{0,\ep}, \ \vc{u}(0, \cdot) =
\vc{u}_{0, \ep}, \ \Div \vc{u}_{0, \ep} = 0.
\]
In this paper, we focus on the \emph{ill-prepared} data, specifically,
\bFormula{ipd1}
n_e(0, \cdot) = N_i + \ep N_{0, \ep}, \ \{ N_{0, \ep} \}_{\ep > 0}
\ \mbox{bounded in} \ L^2 \cap L^\infty (\Omega),
\eF
\bFormula{ipd2}
\vc{u}(0, \cdot) = \vc{u}_{0, \ep}, \ \{ \vc{u}_{0, \ep} \}_{\ep > 0}
\ \mbox{bounded in} \ L^2(\Omega;R^3).
\eF
In particular, the gradient part of the velocity field will develop fast oscillations in the asymptotic limit $\ep \to 0$.

\subsection{Spatial domain}
\label{SD}

In contrast with \cite{AliCheJun}, we consider the physically relevant
(unbounded) domains with boundaries. Similarly to Farwig, Kozono, and Sohr \cite{FAKOSO}, we focus on the class of \emph{uniform}
$C^{3}$ domains of type $(\alpha, \beta, K)$. Specifically,
for each point of $x_0 \in \partial \Omega$, there is a function $h \in C^{3}(R^2)$,
$\| h \|_{C^{3}(R^2)} \leq K$,
and
\[
U_{\alpha, \beta,h} = \{ (y, x_3) \ | h(y) - \beta < x_3 < h(y) + \beta ,\
|y| < \alpha \}
\]
such that, after suitable translation and rotation of the coordinate axes,
$x_0 = [0,0,h(0)]$ and
\[
\Omega \cap U_{\alpha, \beta, h} = \{ (y, x_3) \ | h(y) - \beta < x_3 < h(y) ,\
|y| < \alpha \},
\]
\[
\partial \Omega \cap U_{\alpha, \beta, h} =
\{ (y, x_3) \ |  x_3 = h(y) ,\
|y| < \alpha \}.
\]

Additional hypotheses imposed on the class of domains are stronger for the 
inviscid limit so we consider the two cases separately.

\subsubsection{Hypotheses in the case of constant viscosity}
\label{domvisc}

Since our method leans essentially on dispersion of 
acoustic waves, we suppose that

\medskip

$\bullet$ the \emph{point spectrum} of the Neumann Laplacian $\Delta_N$ in
$L^2(\Omega)$ is empty,

\medskip

\noindent in particular, the domain $\Omega$ must be unbounded.
Although the absence of eigenvalues for the Neumann Laplacian represents, in
general, a delicate and highly unstable problem (see Davies and Parnovski \cite{DavPar}),
there are numerous examples of such domains - the whole space $R^3$,
the half-space, exterior domains, unbounded strips, tube-like domains and waveguides, see
D'Ancona and Racke \cite{DanRac}.

\subsubsection{Hypotheses in the case of inviscid limit}

\label{dominvisc}

The absence of eigenvalues for the Neumann Laplacian is apparently not sufficient to carry over the inviscid limit. We need stronger dispersion provided by the so-called $L^1-L^\infty$ estimates well-known for the acoustic equation in $R^3$, cf. Section \ref{decosc} below. More specifically, we focus on the class of physically relevant domains represented by \emph{infinite waveguides} in the spirit of D'Ancona and Racke \cite{DanRac}. We suppose that 
\bFormula{waveg1}
\Omega \subset R^3, \ \Omega = B \times R^L ,\ L=1,2,3,
\eF
where 
\bFormula{waveg2}
B \subset R^{3 - L} \ \mbox{is a smooth bounded domain for}\ L=1,2, 
\ \Omega = R^3 \ \mbox{for}\ L = 3.
\eF

Obviously, the domains satisfying (\ref{waveg1}), (\ref{waveg2}) belong to the class of uniform $C^3$ domains of type $(\alpha, \beta, K)$, and the point spectrum of the Neumann Laplacian is empty. A peculiar feature of the present problem is that propagation of \emph{acoustic waves} is governed by a wave equation of \emph{Klein-Gordon} type (see Section \ref{a}), where dispersion is enhanced by the presence of ``damping''. In particular, we recover the $L^1-L^\infty$ estimates even in the case of infinite tubes 
($L = 1$) under the Neumann boundary conditions, see Section \ref{decosc}.   

\subsection{Asymptotic limit}

By analogy with the low Mach number limits, we expect the limit velocity
to satisfy the incompressible Navier-Stokes system with a Brinkman type
damping if $\mu_\ep = {\rm const} > 0$, and the Euler system (\ref{euler1} -
\ref{euler3}) in the inviscid limit $\mu_\ep \to 0$.

In comparison with the low Mach number limit, the main difficulty here is the presence of the extra term
\[
\frac{1}{\ep^2} n_e \Grad \Phi = \frac{ n_e - N_i }{\ep}
\Grad \Delta^{-1}_N \left[ \frac{ n_e - N_i }{\ep} \right] +
\frac{N_i}{\ep^2} \Grad \Phi
\]
in the momentum equation (\ref{i2}). While the gradient component
$\frac{N_i}{\ep^2} \Grad \Phi$ can be easily incorporated into the pressure
in the limit system, the quantity
\[
\frac{ n_e - N_i }{\ep}
\Grad \Delta^{-1}_N \left[ \frac{ n_e - N_i }{\ep} \right]
\]
should ``disappear'' in the course of the limit process $\ep \to 0$.
To achieve this, the dispersive estimates based on the celebrated RAGE theorem will be used.

Another difficulty lies in the fact that the quantity
\[
\frac{ n_e - N_i }{\ep}
\Grad \Delta^{-1}_N \left[ \frac{ n_e - N_i }{\ep} \right]
\]
is (known to be) only locally integrable; for global analysis, it must be written in the form
\[
\frac{ n_e - N_i }{\ep}
\Grad \Delta^{-1}_N \left[ \frac{ n_e - N_i }{\ep} \right] =
\frac{1}{\ep^2} \left( \Div \left( \Grad \Phi \otimes \Grad \Phi \right) -
\frac{1}{2} \Grad | \Grad \Phi |^2 \right),
\]
meaning as an element of the dual space $W^{-1,1}$.

Last but not least, we point out that the analysis of the \emph{inviscid} limit leans heavily on the fact that the propagation of acoustic waves is governed by the Klein-Gordon wave equation, yielding effective dispersion on the waveguide like domains specified in Section \ref{dominvisc}. 

The paper is organized as follows. In Section \ref{sw}, we introduce the concept of \emph{suitable weak solution} to system (\ref{i1} - \ref{i7})
that proved to be very convenient for studying the inviscid limits, cf.
\cite{FNP5}. Section \ref{m} contains the main results. In Section \ref{ub},
we summarize the uniform bounds independent of the scaling parameter $\ep$.
Section \ref{a} is devoted to the acoustic equation and the resulting dispersive estimates. Finally, in Section \ref{c}, we show convergence toward the incompressible Navier-Stokes system in the case
of non-degenerate viscosity, while Section \ref{z} completes the proof of the inviscid limit.

\section{Suitable weak solutions}

\label{sw}

Motivated by the genera theory developed by Ruggeri and Trovato \cite{RUTR}, we assume that the electron pressure
$p$ satisfies
\bFormula{pres}
p \in C[0,\infty) \cap C^2(0, \infty),\ p(0)= 0,\
p'(n) > 0 \ \mbox{for all}\ n > 0, \ \lim_{n \to \infty}
\frac{p'(n)}{n^{2/3}} = p_\infty > 0.
\eF

Next, we introduce the standard \emph{Helmholtz decomposition} of a vector
field $\vc{v}$,
\[
\vc{v} = \vc{H}[ \vc{v} ] + \vc{H}^\perp [\vc{v}],
\]
with
\[
\vc{H}^\perp = \Grad \Psi, \ \Delta \Psi = \Div \vc{v} \ \mbox{in}\
\Omega, \ (\Grad \Psi - \vc{v}) \cdot \vc{n}|_{\partial \Omega} = 0.
\]
As shown by Farwig, Kozono and Sohr \cite{FAKOSO}, the linear operator
\[
\vc{H} \ \mbox{is bounded in} \ L^2 \cap L^q (\Omega; R^3)
\ \mbox{for}\ 2 < q < \infty, \ \mbox{and in}\
L^2 + L^q (\Omega) \ \mbox{for} \ 1 < q \leq 2
\]
as soon as $\Omega$ is a $C^2-$domain of type $(\alpha, \beta, K)$ introduced in Section \ref{SD}. Moreover, the norm of $\vc{H}$ in the aforementioned spaces depends solely on the parameters $(\alpha, \beta, K)$. As a matter of fact, the domains considered in the present paper belong to the higher regularity class $C^3$ for several technical reasons that will become clear in the course of the proof of the main results.

Following \cite{FENOSU} we say that a trio $n_e$, $\vu$, $\Phi$ is a \emph{suitable weak solution} to system (\ref{i1} - \ref{i7}) in
$(0,T) \times \Omega$, supplement with the initial conditions $n_e(0, \cdot) =
n_0$, $\vu(0, \cdot) = \vu_0$ if:

\medskip \hrule \medskip

\begin{itemize}
\item
the functions $n_e$, $\vu$, $\Phi$ belong to the regularity class
\[
n_e \geq 0, \ n_e - N_i \in L^\infty(0,T;
L^{5/3} +L^2(\Omega)),
\]
\[
p(n_e) \in L^1(0,T; L^1_{\rm loc} (\Omega)),
\]
\[
\vu \in L^2(0,T; W^{1,2}(\Omega; R^3)),\;\vc u\cdot\vc
n|_{\partial\Omega}=0,
\]
\[
\Grad \Phi \in L^\infty(0,T; L^2(\Omega; R^3));
\]

\item equation of continuity (\ref{i1}) is satisfied in the sense of renormalized solutions
(see DiPerna and P.-L.Lions \cite{DL}),
\bFormula{w1} \int_0^T
\intO{ \left( \left( b(n_e) + n_e \right) \partial_t \varphi +
\left( b(n_e) + n_e \right) \vu \cdot \Grad \varphi + \Big( b(n_e) -
b'(n_e) n_e \Big) \Div \vu \varphi \right) } \ \dt \eF
\[
= - \intO{ \left( b(n_0) + n_0 ) \right) \varphi (0,\cdot) }
\]
for any test function $\varphi \in \DC([0,T) \times \Ov{\Omega})$
and any $b \in C^{\infty}[0,\infty)$, $b' \in \DC[0,\infty)$;

\item momentum equation (\ref{i2}), together with the slip boundary
condition (\ref{i5}), is satisfied in a weak sense,
\bFormula{w2}
\int_0^T \intO{ \Big( m_e n_e \vu \cdot \partial_t \varphi + m_e n_e (\vu
\otimes \vu) : \Grad \varphi + p(n_e) \Div \varphi \Big) } \ \dt \eF
\[
= \int_0^T\intO{\Big(\tn{S}(\Grad \vu) :
\Grad \varphi - n_e \Grad \Phi \cdot \varphi  +
\frac{m_e n_e}{\tau} \vu \cdot \varphi \Big)} \ \dt - \intO{ n_0 \vc{u}_0 \cdot \varphi (0,
\cdot) }
\]
for any test function $\varphi \in \DC([0,T) \times \Ov{\Omega}; R^3)$,
$\varphi\cdot\vc n|_{[0,T)\times\partial\Omega}=0$;
\item the electric potential $\Phi$ is given by formula
\bFormula{w2a}
\Grad \Phi (s, \cdot) = \Grad \Phi_{0} - \int_0^s \vc{H}^\perp [
n_e \vu ] \ {\rm d}t,
\eF
where
\bFormula{w2aa}
\Delta \Phi_{0} = n_0 - N_i \ \mbox{in} \ \Omega, \ \Grad \Phi_0 \cdot
\vc{n} |_{\partial \Omega} = 0;
\eF
\item the \emph{relative entropy inequality}
\bFormula{w3}
\intO{ \left( \frac{1}{2} m_e n_e |\vu - \vc{U}|^2 + E(n_e, r)
 + \frac{1}{2} |\Grad \Phi |^2 \right) (s, \cdot)  }
\eF
\[
+ \int_0^s \intO{ \left[ \tn{S} (\Grad \vu) -
\tn{S}(\Grad \vc{U}) \right] : \Grad (\vc{u} - \vc{U}) } \ \dt
+
\int_0^s \intO{ \frac{m_e}{\tau} n_e | \vu - \vc{U} |^2 } \ \dt
\]
\[
\leq \intO{ \left( \frac{1}{2} m_e n_0 |\vu_0 - \vc{U}(0, \cdot)|^2 +
E(n_0, r(0,\cdot)) + \frac{1}{2} |\Grad \Phi_0 |^2  \right) } + \int_0^s {\cal R} \left( n_e,
\vu, r, \vc{U} \right) \ \dt
\]
holds for a.a. $s \in [0,T]$ and all test functions $r$, $\vc{U}$
such that
\[
r - N_i \in \DC ( [0,T] \times \Ov{\Omega} ),\ r > 0, \ \vc{U}
\in \DC([0,T] \times \Ov{\Omega};R^3), \ \vc{U} \cdot
\vc{n}|_{\partial \Omega} = 0,
\]
where
\bFormula{w4}
{\cal R}\left( n_e, \vu, r, \vc{U} \right) \equiv
\intO{  m_e n_e \Big( \partial_t \vc{U} + \vu \Grad \vc{U} \Big) \cdot
(\vc{U} - \vu )}
\eF
\[
+ \intO{\tn{S}(\Grad \vc{U}):\Grad (\vc U- \vc{u})  } +
\intO{ \frac{m_e}{\tau} n_e \vc{U} \cdot (\vc{U} - \vu ) }
- \intO{n_e \Grad \Phi \cdot \vc{U}}
\]
\[
+ \intO{ \left( (r - n_e) \partial_t P(r) + \Grad P(r) \cdot
\left( r \vc{U} - n_e \vu \right) - \Div \vc{U} \Big( n_e \Big(
P(n_e) - P(r) \Big) - E (n_e, r) \Big)\right) },
\]
with
\[
E(n_e, r) \equiv H(n_e) - H'(r) (n_e - r) - H(r).
\]
and
\[
P \equiv H',\; H(n) \equiv n \int_1^n \frac{p(s)}{s^2}{\rm d}s.
\]

\end{itemize}

\medskip \hrule \medskip

It can be deduced from (\ref{w2a}) that
\bFormula{w2aaa}
\intO{ \Grad \Phi (t, \cdot) \Grad \varphi } = \intO{ (N_i - n_e) (t, \cdot) \varphi }
\ \mbox{for any} \ t \in (0,T) \ \mbox{and all test functions}\
\varphi \in \DC (\Ov{\Omega}),
\eF
meaning $\Phi$ is a (strong) solution of Poisson equation (\ref{i3}).
In particular, by virtue of the standard (local) elliptic regularity,
\bFormula{w2A}
\Phi(t, \cdot) \in W^{2,5/3}(K) \ \mbox{for any compact} \ K \subset \Ov{\Omega}.
\eF

The \emph{existence} of global-in-time suitable weak solutions to the
compressible Navier-Stokes system in a bounded spatial domain and the no-slip
boundary conditions was proved in \cite[Theorem 3.1]{FENOSU} with the help of
the approximation scheme introduced in \cite{FNP}. Adaptation of the method to the present problem requires only straightforward modifications. The main advantage of working directly with suitable weak solutions is that the relative entropy inequality (\ref{w3}) already implicitly includes the stability estimates necessary to perform the inviscid limit.

\section{Main results}

\label{m}

We start by introducing the scaled system. To simplify notation, we set
$m_e = \ep^2$, $n_e = n_\ep$, $\vu = \vue$, $\Phi = \Phi_\ep$, and
$N_i = \Ov{n}$ - a positive constant. Accordingly, system of equations
(\ref{i1} - \ref{i3}) reads
\bFormula{si1}
\partial_t \nep + \Div (\nep \vue)=0
\eF
\bFormula{si2}
\partial_t ( \nep \vue) + \Div( \nep \vue \otimes \vue )  +
\frac{1}{\ep^2} \Grad p(\nep) = \Div \tn{S}_\ep (\Grad \vue) + \frac{1}{\ep^2} \nep \Grad \Phi_\ep - \frac{1}{\tau} \nep \vue,
\eF
\bFormula{si3}
\Delta \Phi_\ep = \nep - \Ov{n},
\eF
with the viscous stress
\bFormula{si4}
\tn{S}_\ep(\Grad \vue) = {\mu_\ep} \Big( \Grad \vue + \Grad^t \vue -\frac{2}{3} \Div \vue \tn{I} \Big), \ \mu_\ep  > 0.
\eF
System (\ref{si1} - \ref{si3}) is supplemented with the boundary conditions
\bFormula{si5}
\vue \cdot \vc{n}|_{\partial \Omega} = 0,\
[ \tn{S}_\ep (\Grad \vue) \vc{n} ]  \times \vc{n}|_{\partial \Omega} = 0,
\eF
\bFormula{si6}
\Grad \Phi_\ep \cdot \vc{n}|_{\partial \Omega} = 0,
\eF
and
\bFormula{si7}
\vue(x) \to 0 , \ \Phi_\ep(x) \to 0  \ \mbox{as}\ |x| \to \infty.
\eF

\subsection{Ill-prepared initial data}

Taking
\[
r \equiv \Ov{n}, \ \vc{U} \equiv 0
\]
as test functions
in the relative entropy inequality (\ref{w3}) we obtain

 \bFormula{ub1}
\intO{ \left( \frac{1}{2} \nep |\vue |^2 + \left[ \frac{
H(\nep) - H'(\Ov{n})(\nep - \Ov{n}) - H(\Ov{n})}{\ep^2} \right]
 + \frac{1}{2 \ep^2} |\Grad \Phi_\ep |^2 \right) (s, \cdot)  }
\eF
\[
+ \int_0^s \intO{ \left( {\mu_\ep}
\left| \Grad \vue + \Grad^t \vue - \frac{2}{3} \Div \vue \tn{I} \right|^2   +
\frac{1}{\tau} \nep |\vue|^2  \right) } \ \dt
\]
\[
\leq \intO{ \left( \frac{1}{2} n_{0,\ep} |\vu_{0,\ep}|^2 +
\frac{1}{\ep^2} E(n_{0,\ep}, \Ov{n}) + \frac{1}{2 \ep^2} |\Grad \Phi_{0,\ep} |^2  \right) },
\]
where
\bFormula{inda}
n_\ep = n_{0,\ep}, \ \vue(0, \cdot) = \vu_{0, \ep}
\eF
and
\bFormula{inda2}
\Phi_{0, \ep} = \Delta^{-1}_N [ n_{0,\ep} - \Ov{n}]
\eF
are the initial data.

Consequently, the initial data must be chosen in such a way that the expression on the right-hand side of (\ref{ub1}) remains bounded for $\ep \to 0$. Accordingly, we suppose that
\bFormula{id1}
n_\ep (0, \cdot) = n_{0,\ep} = \Ov{n} + \ep N_{0,\ep}, \
\{ N_{0,\ep} \}_{\ep > 0} \ \mbox{bounded in} \ L^2 \cap L^\infty (\Omega),
\eF
\bFormula{id2}
\vue(0, \cdot) = \vu_{0, \ep}, \ \{ \vu_{0, \ep} \}_{\ep > 0} \
\mbox{bounded in}\ L^2(\Omega;R^3).
\eF
Moreover, the functions $N_{0,\ep}$ must be taken so that
$\Phi_{0,\ep} = \ep \Delta^{-1}_N [ N_{0,\ep} ]$ satisfy
\bFormula{id3}
\{ \Grad \Phi_{0,\ep} \}_{\ep > 0} \ \mbox{bounded in} \ L^2(\Omega;R^3).
\eF

\subsection{Asymptotic limit for positive viscosity coefficients}

Our first result concerns the asymptotic limit in the case $\mu_\ep =
\mu > 0$.

\bTheorem{m1}
Let $\Omega \subset R^3$ be an (unbounded) $C^3-$domain of type
$(\alpha, \beta,K)$ specified in Section \ref{SD} and such that the point
spectrum of the Neumann Laplacian $\Delta_N$ in $L^2(\Omega)$ is empty.
Suppose that the viscosity coefficient $\mu_\ep = \mu > 0$ is independent of
$\ep$ and that the pressure $p$ satisfies (\ref{pres}).
Let $\{ n_\ep, \vue, \Phi_\ep \}_{\ep > 0}$ be a sequence of suitable weak solutions to the scaled system (\ref{si1} - \ref{si7}), emanating from the initial data satisfying (\ref{id1} - \ref{id3}).

Then
\bFormula{m1}
{\rm ess} \sup_{t \in (0,T)} \| \nep - \Ov{n} \|_{L^2 + L^{5/3}(\Omega)} \leq \ep c,
\eF
and, at least for a suitable subsequence,
\[
\vu_{0,\ep} \to \vc{U}_0 \ \mbox{weakly in}\ L^2(\Omega;R^3),
\]
\bFormula{m2}
\vue \to \vc{U} \ \mbox{weakly in} \ L^2(0,T; W^{1,2}(\Omega;R^3))
\ \mbox{and (strongly) in} \ L^2((0,T) \times K;R^3)
\eF
for any compact $K \subset \Omega$, where $\vc{U}$ is a weak solution to the incompressible (damped) Navier-Stokes system in $(0,T) \times \Omega$,
\bFormula{m3}
\Div \vc{U} = 0,
\eF
\bFormula{m4}
\Ov{n} \Big( \partial_t \vc{U} + \vc{U} \cdot \Grad \vc{U} \Big) + \Grad \Pi = \Div \tn{S}(\Grad \vc{U}) - \frac{\Ov{n}}{\tau} \vc{U},
\eF
with
\bFormula{m5}
\vc{U} \cdot \vc{n} |_{\partial \Omega} = 0, \ [\tn{S}(\Grad \vc{U}) \vc{n}]
\times \vc{n}|_{\partial \Omega} = 0,
\eF
and
\bFormula{m6}
\vc{U}(0, \cdot) = \vc{H}[\vc{U}_0].
\eF
\eT

\medskip

{\bf Remark \ref{m}.1}  {\it Momentum equation (\ref{m4}), together with
the slip boundary conditions (\ref{m5}) and the initial condition (\ref{m6}), are understood in the weak sense, specifically, the integral identity
\bFormula{m7}
\int_0^T \intO{ \Ov{n} \Big( \vc{U} \cdot \partial_t \varphi +
(\vc{U} \otimes \vc{U}) : \Grad \varphi \Big) } \ \dt =
\int_0^T \intO{ \Big( \tn{S}(\Grad \vc{U}): \Grad \varphi +
\frac{\Ov{n}}{\tau} \vc{U} \cdot \varphi \Big) } \ \dt -
\intO{ \vc{U}_0 \cdot \varphi (0,\cdot) }
\eF
for any test function $\varphi \in \DC([0,T) \times \Ov{\Omega}; R^3)$,
$\Div \varphi = 0$, $\varphi \cdot \vc{n}|_{\partial \Omega} = 0$.
}

\subsection{Inviscid limit}

Our second result concerns the case of vanishing viscosity coefficient
$\mu_\ep \searrow 0$. In this case, the limit velocity field is expected to satisfy the incompressible Euler system (\ref{euler1} - \ref{euler3}).
As is well-known, this system
possesses a local-in-time solution
\bFormula{EU1}
\vc{v} \in C([0,T_{\rm max}); W^{k,2}(\Omega)),\ \Grad \Pi \in
C([0,T_{\rm max}); W^{k-1,2}(\Omega))
\eF
provided
\[
\vc{v}(0, \cdot) = \vc{v}_0 \in W^{k,2}(\Omega), \ k > 5/2,\ \Div
\vc{v}_0 = 0,\ \vc{v}_0 \cdot \vc{n}|_{\partial \Omega} = 0,
\]
and provided $\Omega = R^3$, $\Omega$ is a half-space, or $\Omega$ is an
(exterior) domain with compact boundary.
The life-span $T_{\rm max}$ depends solely on $\| \vc{v}_0
\|_{W^{k,2}(\Omega;R^3)}$, see Alazard \cite{AL2}, Isozaki
\cite{Isoz}, Secchi \cite{Secc}, among others. As a matter of fact,
the damping term $\frac{1}{\tau} \vc{v}$ in (\ref{euler2}) may extend the life-span of regular solutions, in particular if the initial data are small in comparison with $1/\tau$. In a very interesting
recent paper, Chae \cite{Chae} showed that a smooth solution of (\ref{euler1} - \ref{euler3}) exists \emph{globally} in time provided $\Omega = R^3$, and
$\tau < T^E_{\rm max}$, where $T^E_{\rm max}$ is the life span of the regular solution of the undamped Euler system emanating from the same initial data.

\bTheorem{m2}
In addition to hypotheses of Theorem \ref{Tm1},
assume that $\Omega \subset R^3$ is an infinite waveguide specified in Section 
\ref{dominvisc}. Moreover, we suppose that
\[
\mu_\ep \searrow 0,
\]
and that the initial data satisfy
\[
N_{0,\ep} \to N_0 \ \mbox{(strongly) in}\ L^2(\Omega),\
\vu_{0,\ep} \to \vc{u}_0 \ \mbox{(strongly) in} \ L^2(\Omega;R^3))
\]
as $\ep \to 0$, where
\[
\vc{v}_0 = \vc{H}[\vc{u}_0] \in \ W^{k,2}(\Omega;R^3) , \ k > 5/2.
\]
Moreover, suppose that the damped Euler system (\ref{euler1}-
\ref{euler3}), with the initial datum $\vc{v}_0$, possesses a regular solution $\vc{v}$ defined on a time interval $[0, T_{\rm max})$ satisfying
(\ref{EU1}).

Then
\bFormula{mm1}
{\rm ess} \sup_{t \in (0,T)} \| \nep - \Ov{n} \|_{L^2 + L^{5/3}(\Omega)} \leq \ep c,
\eF
and
\[
{\rm ess} \sup_{t \in (0,T_{\rm loc})} \left\| \vc{H}\left[
\sqrt{\frac{\nep}{{\Ov{n}}}} \vue \right] - \vc{v}
 \right\|_{L^2 (\Omega;R^3)} \to 0,
\]
\[
\sqrt{\frac{\nep}{{\Ov{n}}}} \vue \to \vc{v} \ \mbox{in}\
L^{q} (0,T_{\rm loc}; L^2(K; R^3)) \ \mbox{for any compact}\
K \subset \Omega, \ \ 1 \leq q < \infty,
\]
for any $T_{\rm loc} < T_{\rm max}$, $T_{\rm loc} \leq T$.
\eT

\medskip

{\bf Remark \ref{m}.2}
{\it
The proof of Theorem \ref{Tm2} leans essentially on the $L^1 - L^\infty$ bounds for acoustic waves established in Section \ref{decosc}. Thus the conclusion of Theorem \ref{Tm2} remains valid as soon as 
these bounds are available. Note that Isozaki \cite{Isoz} established similar estimates on \emph{exterior} domains in $R^3$. 
}

\medskip

The rest of the paper is devoted to the proof of Theorems \ref{Tm1}, \ref{Tm2}.

\section{Uniform bounds}
\label{ub}

For the ill-prepared initial data,
all desired uniform bounds follow from the energy inequality (\ref{ub1}).
Introducing the \emph{essential} and \emph{residual} parts of a quantity
$h$,
\[
h = [h]_{\rm ess} + [h]_{\rm res}
\]
\[
[h]_{\rm ess} = \chi (n_\ep) h, \ [h]_{\rm res} = h - [h_{\rm ess}],
\]
where
\[
\chi \in \DC(0,\infty), \ 1 \leq \chi \leq 1 , \ \chi(n) = 1
\ \mbox{for}\ n \ \mbox{belonging to an open neighborhood of}\
\Ov{n} \ \mbox{in}\ (0, \infty),
\]
we get the following list of estimates
\bFormula{ub2}
{\rm ess} \sup_{t \in (0,T)} \| \sqrt{\nep} \vue (t, \cdot) \|_{L^2(\Omega;R^3)} \leq c,
\eF
\bFormula{ub3}
{\rm ess} \sup_{t \in (0,T)} \| \Grad \Phi_{\ep} (t, \cdot) \|_{L^2(\Omega;R^3)} \leq \ep c,
\eF
\bFormula{ub4}
{\rm ess} \sup_{t \in (0,T)} \left\| \left[ \frac{ \nep - \Ov{n} }{\ep}
\right]_{\rm ess} (t, \cdot) \right\|_{L^2(\Omega)} \leq c ,\
{\rm ess} \sup_{t \in (0,T)} \intO{ [\nep]^{5/3}_{\rm res} (t,\cdot) } \leq \ep^2 c,
\eF
\bFormula{ub6}
{\rm ess} \sup_{t \in (0,T)} \intO{ [1]_{\rm res} (t, \cdot) } \leq \ep^2 c,
\eF
and
\bFormula{ub7}
\mu_\ep \int_0^T \intO{ \left| \Grad \vue + \Grad \vue^t - \frac{2}{3}
\Div \vue \tn{I} \right|^2 } \ \dt \leq c,
\eF
where all generic constants are independent of $\ep$.

Estimates (\ref{ub4} - \ref{ub7}) can be combined to deduce a bound on the velocity field in the Sobolev space $L^2(0,T; W^{1,2}(\Omega;R^3))$ that is relevant in the proof of Theorem \ref{Tm1}. To this end, we report the following version of \emph{Korn's inequality} that may be of independent interest.

\bProposition{ub}
Let $\Omega \subset R^3$ be a $C^2$-uniform domain of type $(\alpha, \beta, K)$ introduced in Section \ref{SD}.

Then there exists $\delta > 0$, depending solely on the parameters $(\alpha, \beta, K)$, such that
\bFormula{korn}
\| \vc{w} \|_{W^{1,2}(\Omega;R^3)}^2 \leq
c (\alpha, \beta, \delta, K) \left( \left\| \Grad \vc{w} + \Grad^t \vc{w} -
\frac{2}{3} \Div \vc{w} \tn{I} \right\|^2_{L^2(\Omega; R^{3 \times 3})} +
\int_{\Omega \setminus V} |\vc{w}|^2 \ \dx \right)
\eF
for any measurable set $V$, $|V| < \delta$, and for all $\vc{w} \in W^{1,2}
(\Omega;R^3)$.
\eP

\bProof

In view of the standard decomposition technique and partition of unity, it is enough to show (\ref{korn}) on each set
\[
U^-_{\alpha, \beta, K} = \{ (y, x_3) \ | h(y) - \beta < x_3 < h(y) ,\
|y| < \alpha \}.
\]

Revoking the result \cite[Proposition 4.1] {BuFe1}, we have
\bFormula{korn1}
\| \vc{w} \|_{W^{1,2}(U^-;R^3)}^2 \leq
c (\alpha, \beta, K) \left( \left\| \Grad \vc{w} + \Grad^t \vc{w} -
\frac{2}{3} \Div \vc{w} \tn{I} \right\|^2_{L^2(U^-; R^{3 \times 3})} +
\int_{U^-} |\vc{w}|^2 \ \dx \right).
\eF
As a matter of fact, the constant $c$ in (\ref{korn1}) depends only
the Lipschitz constant of the function $h$ and width of $U^-$ given in terms of $\alpha$, $\beta$.

Furthermore, we have
\[
| U^-_{\alpha, \beta, K} | \geq 2 \delta  > 0
\]
for a certain $\delta (\alpha, \beta, K) > 0$. In particular,
\[
|U^- \setminus V| > \delta \ \mbox{for any measurable set}\ V, \
|V| < \delta.
\]

Now, arguing by contradiction, we construct sequences
\[
\{ h_n \}_{n=1}^\infty, \ \| h_n \|_{C^2(R^2)} \leq K,\
U^-_n = \Big\{ (y, x_3) \ \Big| \ h_n (y) - \beta < x_3 < h_n(y), \ |y| < \alpha \Big\},
\]
\[
h_n \to h \ \mbox{in} \ C^1( \{ |y| \leq \alpha ),
\]
\[
\{ V_n \}_{n=1}^\infty,\ |V_n| < \delta,
\]
and
\[
\{ \vc{w}_n \}_{n=1}^\infty, \ \| \vc{w}_n \|_{W^{1,2}(U^-_n;R^3)} = 1,
\]
such that
\[
\left( \left\| \Grad \vc{w}_n + \Grad^t \vc{w}_n -
\frac{2}{3} \Div \vc{w}_n \tn{I} \right\|^2_{L^2(U^-_n; R^{3 \times 3})} +
\int_{U^-_n \setminus V_n} |\vc{w}_n|^2 \ \dx \right) \to 0 \ \mbox{as}\
n \to \infty.
\]

Because the domains $U^-_n$ are uniformly Lipschitz, we can extend $\vc{w}_n$
as $\tilde \vc{w}_n$ on the cylinder
\[
U = \{ (y, x_3) \ | \ |x_3| < 2 + \beta + \alpha K,\ |y| < \alpha \}
\]
in such a way that
\[
\| \tilde \vc{w}_n \|_{W^{1,2} (U; R^3)} \leq c(\alpha, \beta, K) \| \vc{w}_n
\|_{W^{1,2}(U^-_n; R^3)}.
\]
Since $W^{1,2}(U ; R^3)$ is compactly embedded into $L^2(U ;R^3)$, we may use (\ref{korn1}) to deduce that
\[
\tilde \vc{w}_n \to \vc{w} \ \mbox{in} \ W^{1,2}(U; R^3),
\]
where
\bFormula{contr1}
\Grad \vc{w} + \Grad^t \vc{w} - \frac{2}{3} \Div \vc{w} \tn{I} = 0 ,\
\vc{w} \not\equiv 0
\ \mbox{in the set} \ \Big\{ (y, x_3) \ \Big| \ h(y) - \beta < x_3 <
h(y) , \ |y| < \alpha \Big\}.
\eF
On the other hand,
\bFormula{contr2}
\int_U 1_{U^-_n \setminus V_n} |\vc{w}_n|^2 \ \dx \to
\int_U \chi |\vc{w}|^2 \ \dx = 0
\eF
where
\[
1_{U^-_n \setminus V_n} \to \chi \ \mbox{weakly-(*) in} \ L^\infty (\Omega),
\ \chi \geq 0, \ \int_{U} \chi \ \dx > 0.
\]
However, relation (\ref{contr1}) implies that $\vc{w}$ is a (nonzero) conformal Killing vector (see Reshetnyak \cite{Resh}) vanishing, by virtue (\ref{contr2}), on a set of positive measure, which is impossible.

\qed

Thus, finally, taking $V$ the ``residual set'',
$V= {\rm supp} [1]_{\rm res}$ we may combine the estimates (\ref{ub1}),
(\ref{ub6}), and (\ref{ub7}) with Proposition \ref{Pub} to conclude that
\bFormula{ub7a}
\mu_\ep \int_0^T  \| \vue \|^2_{W^{1,2}(\Omega;R^3)} \ \dt \leq c.
\eF

\section{Acoustic equation}
\label{a}

As already pointed out, the essential piece of information
necessary to carry out the asymptotic limit is contained in the
oscillatory component of the velocity field responsible for propagation of acoustic waves. Introducing new variables
\[
N_\ep = \frac{\nep - \Ov{n}}{\ep}, \ \vc{V}_\ep = \nep \vue
\]
we can formally rewrite system (\ref{si1}), (\ref{si2}) in the form
\bFormula{a1}
\ep \partial_t N_\ep + \Div \vc{V}_\ep = 0,
\eF

\bFormula{a2}
\ep \partial_t \vc{V}_\ep + p'(\Ov{n}) \Grad N_\ep +
\frac{\Ov{n}}{\ep} \Grad \Phi_\ep
\eF
\[
= \ep \Div \tn{S} (\Grad \vue) - \ep \Div( \nep \vue \otimes \vue)
- \frac{\ep}{\tau} n_\ep \vue
- N_\ep \Grad \Phi_\ep - \frac{1}{\ep} \Grad \Big(
p(\nep) - p'(\Ov{n}) (\nep - \Ov{n}) - p(\Ov{n}) \Big),
\]
supplemented with the boundary condition
\[
\vc{V}_\ep \cdot \vc{n}|_{\partial \Omega} = 0.
\]

System (\ref{a1}), (\ref{a2}) is usually called \emph{acoustic equation},
see Lighthill \cite{LIG}. Its (rigorous) weak formulation reads
\bFormula{wa1}
\int_0^T \intO{ \Big( \ep N_\ep \partial_t \varphi + \vc{V}_\ep
\cdot \Grad \varphi \Big) } \ \dt = - \ep \intO{ N_{0,\ep} \varphi (0, \cdot) }
\eF
for any $\varphi \in \DC([0,T) \times \Ov{\Omega})$, and
\bFormula{wa2}
\int_0^T \intO{ \Big( \ep \vc{V}_\ep \cdot \partial_t \Grad \varphi +
p'(\Ov{n}) N_\ep \Delta \varphi - \Ov{n} N_\ep \varphi \Big) } \ \dt
\eF
\[
= \ep \int_0^T \intO{  \tn{S}_\ep (\Grad \vue) : \Grad^2 \varphi  } \ \dt
- \ep \int_0^T \intO{ \left( \nep \vue \otimes \vue +
\frac{ p(\nep) - p'(\Ov{n})(\nep - \Ov{n}) - p(\Ov{n}) }{\ep^2} \tn{I}
\right) : \Grad^2 \varphi } \ \dt
\]
\[
+ \ep \int_0^T \intO{ \left( \frac{1}{\tau} \nep \vue \cdot \Grad \varphi
+ N_\ep \Grad \left(\frac{\Phi_\ep}{\ep} \right) \cdot \Grad \varphi \right) } \ \dt  - \ep \intO{ n_{0,\ep} \vu_{0,\ep} \cdot \Grad \varphi (0, \cdot) } \]
for any $\varphi \in \DC(\Omega)$, $\Grad \varphi \cdot \vc{n}|_{\partial \Omega} = 0$. Moreover, we rewrite
\[
\intO{ N_\ep \Grad \left( \frac{\Phi_\ep}{\ep} \right) \cdot \Grad \varphi } =
\frac{1}{\ep^2} \intO{ \Delta \Phi_\ep \Grad \Phi_\ep \cdot \Grad \varphi }
\]
\[
= \frac{1}{\ep^2} \intO{ \left( \frac{1}{2} \Grad |\Grad \Phi_\ep |^2 \cdot \Grad \varphi + \Div ( \Grad \Phi_\ep \otimes \Grad \Phi_\ep ) \cdot \Grad \varphi \right)
}
\]
\[
= - \frac{1}{\ep^2} \intO{ \left( \frac{1}{2} |\Grad \Phi_\ep |^2 \cdot \Delta \varphi + ( \Grad \Phi_\ep \otimes \Grad \Phi_\ep ) : \Grad^2 \varphi \right)}.
\]

Furthermore, it follows directly from the uniform bounds established in
(\ref{ub1} - \ref{ub7}) that (\ref{wa2}) can be written as
\bFormula{wa3}
\int_0^T \intO{ \Big( \ep \vc{V}_\ep \cdot \partial_t \Grad \varphi +
p'(\Ov{n}) N_\ep \Delta \varphi - \Ov{n} N_\ep \varphi \Big) } \ \dt
\eF
\[
= \ep \int_0^T \intO{ \Big( \tn{G}^1_\ep : \Grad^2 \varphi +
\tn{G}^2_\ep : \Grad^2 \varphi + \frac{1}{\tau} \nep \vue \cdot \Grad \varphi \Big) } \ \dt - \ep \intO{ n_{0,\ep} \vu_{0,\ep} \cdot \Grad \varphi (0, \cdot) }
\]
for any $\varphi \in \DC(\Omega)$, $\Grad \varphi \cdot \vc{n}|_{\partial \Omega} = 0$, where
\bFormula{aaa1}
\{ \tn{G}^1_\ep \}_{\ep > 0} \ \mbox{is bounded in}\
L^\infty(0,T; L^1(\Omega; R^{3 \times 3})),
\eF
and
\bFormula{aaa2}
\{ \tn{G}^2_\ep \}_{\ep > 0} \ \mbox{is bounded in}\
L^2(0,T; L^2(\Omega; R^{3 \times 3})).
\eF

\subsection{Neumann Laplacian}

At this stage, it is convenient to rewrite the acoustic system (\ref{wa1}),
(\ref{wa3}) in terms of a single self-adjoint operator ${\mathcal A}$
in $L^2(\Omega)$, specifically,
\[
\mathcal{A} = - p'(\Ov{n}) \Delta_N + \Ov{n} \tn{I},
\]
with
\[
{\cal D}(\mathcal{A}) = \Big\{ w \in W^{1,2}(\Omega) \ \Big|
\]
\[
\intO{ \left( p'(\Ov{n}) \Grad w \cdot \Grad \varphi +
\Ov{n} w \varphi \right) } = \intO{ g \varphi }
\ \mbox{for any} \ \varphi \in \DC(\Ov{\Omega})
\ \mbox{for a certain}\ g \in L^2(\Omega) \Big\}.
\]

Given the regularity of the boundary $\partial \Omega$, it can be shown that
\bFormula{NL1}
{\cal D}(\mathcal{A}) = \Big\{ w \in W^{2,2} (\Omega) \ \Big| \
\Grad w \cdot \vc{n} |_{\partial \Omega} = 0 \Big\},\
\mathcal{A}[w] = - p'(\Ov{n}) \Delta w + \Ov{n} w .
\eF
Furthermore, since $\Omega$ is of uniform $C^3$-class, the classical elliptic theory yields
\bFormula{NL2}
{\cal D}(\mathcal{A}^2) \subset C^{2 + \nu} \cap W^{2, \infty} (\Ov{\Omega})
\eF
for a certain $\nu > 0$. We remark that all we need is only uniform $C^{2 + \nu}-$regularity of the boundary instead of $C^3$.

\subsection{Acoustic equation - abstract formulation}

In view of (\ref{NL1}), (\ref{NL2}), and the uniform bounds established in
(\ref{aaa1}), (\ref{aaa2}), the acoustic equation (\ref{wa1}), (\ref{wa3}) can be written in a concise form:
\bFormula{Ab1}
\int_0^T \intO{ \Big( \ep N_\ep \partial_t \varphi + \vc{V}_\ep
\cdot \Grad \varphi \Big) } \ \dt = - \ep \intO{ N_{0,\ep} \varphi (0, \cdot) }
\eF
for any $\varphi \in \DC([0,T) \times \Ov{\Omega})$, and
\bFormula{Ab2}
\int_0^T \intO{ \Big( \ep \vc{V}_\ep \cdot \partial_t \Grad \varphi
- N_\ep \mathcal{A}[\varphi] \Big) } \ \dt = \ep \int_0^T \intO{
F_\ep \mathcal{A}^2 [ \varphi ] }\ \dt - \ep \intO{ n_{0,\ep}
\vu_{0,\ep} \cdot \Grad \varphi (0, \cdot) }
\eF
for any $\varphi \in C^1([0,T); {\cal D}(\mathcal{A}^2))$,
with
\bFormula{Ab3}
\{ F_\ep \}_{\ep > 0}\ \mbox{bounded in}\ L^2(0,T; L^2(\Omega)).
\eF

Thus, using the standard variation-of-constants formula, we obtain
\bFormula{vcf1}
N_{\ep} = \frac{1}{2} \left(
\exp \left( {\rm i} \sqrt{ \mathcal{A} } \frac{t}{\ep} \right) \left[ N_{0, \ep} + \frac{{\rm i}}{\sqrt{ \mathcal{A} }} Z_{0,\ep} \right]
+ \exp \left( - {\rm i} \sqrt{ \mathcal{A} } \frac{t}{\ep} \right) \left[ N_{0, \ep} - \frac{{\rm i}}{\sqrt{ \mathcal{A} }} Z_{0,\ep} \right] \right)
\eF
\[
+ \frac{1}{2} \int_0^t \left(
\exp \left( {\rm i} \sqrt{ \mathcal{A} } \frac{t-s}{\ep} \right) -
\exp \left( - {\rm i} \sqrt{ \mathcal{A} } \frac{t - s}{\ep} \right) \right) \left[  \frac{{\rm i}}{\sqrt{ \mathcal{A} }} \mathcal{A}^2 [F_\ep]  \right]  \ {\rm d}s,
\]
\bFormula{vcf2}
Z_{\ep} = \frac{1}{2} \left(
\exp \left( {\rm i} \sqrt{ \mathcal{A} } \frac{t}{\ep} \right) \left[ Z_{0, \ep} - {{\rm i}}{\sqrt{ \mathcal{A} }} [N_{0,\ep}] \right]
+ \exp \left( - {\rm i} \sqrt{ \mathcal{A} } \frac{t}{\ep} \right) \left[ Z_{0, \ep} + {{\rm i}}{\sqrt{ \mathcal{A} }} [N_{0,\ep}] \right] \right)
\eF
\[
+ \frac{1}{2} \int_0^t \left(
\exp \left( {\rm i} \sqrt{ \mathcal{A} } \frac{t-s}{\ep} \right) +
\exp \left( - {\rm i} \sqrt{ \mathcal{A} } \frac{t - s}{\ep} \right) \right) \left[  \mathcal{A}^2 [F_\ep]  \right]  \ {\rm d}s,
\]
where $Z_\ep$ is interpreted as
\bFormula{Z1}
\intO{ Z_\ep G(\mathcal{A})[ \varphi ] } = - \intO{ \vc{V}_\ep
\cdot \Grad (G(\mathcal{A})[ \varphi ] ) } \ \mbox{for any}\
G \in \DC (\Ov{n}, \infty ), \ \varphi \in \DC(\Omega),
\eF
in particular,
\bFormula{Z2}
\intO{ Z_{0, \ep}  G(\mathcal{A})[ \varphi ] } = - \intO{ n_{0,\ep} \vu_{0, \ep}
\cdot \Grad (G(\mathcal{A})[ \varphi ] )} \ \mbox{for any}\
G \in \DC ( \Ov{n}, \infty  ), \ \varphi \in \DC(\Omega).
\eF
Note that the spectrum of the operator $\mathcal{A}$ is the half-line
$[\Ov{n}, \infty)$.

\subsection{Application of RAGE theorem}
\label{RAJ}

With the explicit formulas (\ref{vcf1}), (\ref{vcf2}) at hand, we are ready to show local energy decay for $N_\ep$ and the acoustic waves represented by the gradient component $\vc{H}^\perp [\vc{V}_\ep]$. To this end,
we employ the following version of the celebrated
RAGE theorem, see Cycon et al. \cite[Theorem 5.8]{CyFrKiSi}:

\bTheorem{rage}
Let $H$ be a Hilbert space, ${A}: {\cal D}({A}) \subset H \to H$ a
self-adjoint operator, $C: H \to H$ a compact operator, and $P_c$
the orthogonal projection onto the space of continuity $H_c$ of $A$,
specifically,
\[
H = H_c \oplus {\rm cl}_H \Big\{ {\rm span} \{ w \in H \ | \ w \ \mbox{an eigenvector of} \ A \} \Big\}.
\]

Then
\bFormula{rage1}
\left\| \frac{1}{\tau} \int_0^\tau \exp(-{\rm i} tA ) C P_c \exp(
{\rm i} tA ) \ \dt \right\|_{{\cal L}(H)} \to 0 \ \mbox{as}\ \tau
\to \infty.
\eF
\eT

We apply Theorem \ref{Trage} to $H= L^2(\Omega)$, $A = - \sqrt{{\mathcal A}}$, $C = \chi^2 G(\mathcal{A})$, with $\chi \in \DC(\Omega)$, $\chi \geq 0$. In accordance with hypotheses of Theorem \ref{Tm1}, the point spectrum
of $\mathcal{A}$ is empty, and we deduce that
\bFormula{rage2}
\int_0^T \left\| \chi G(\mathcal{A})  \exp \left( {\rm i} \sqrt{ \mathcal{A} } \frac{t}{\ep} \right) [X] \right\|^2_{L^2(\Omega)} \ \dt \leq \omega(\ep)
\| X \|_{L^2(\Omega)}^2,
\eF
where $\omega(\ep) \to 0$ as $\ep \to 0$. In particular, going back to
(\ref{vcf1}), (\ref{vcf2}) we may infer that
\[
\left\{ t \mapsto \intO{ N_\ep (t, \cdot) G( \mathcal{A} )[\varphi] }
\right\} \to 0 \ \mbox{in}\ L^2(0,T),
\]
and, similarly,
\[
\left\{ t \mapsto \intO{ \vc{V}_\ep \cdot \Grad ( G( \mathcal{A} )[\varphi])}
\right\} \to 0 \ \mbox{in}\ L^2(0,T)
\]
as $\ep \to \infty$ for any $G \in \DC(\Ov{n}, \infty)$, $\varphi \in \DC(\Omega)$. Thus, by means of density argument,
\bFormula{convv1}
N_\ep \to 0 \ \mbox{in}\ L^q(0,T; (L^2 + L^{5/3})_{\rm weak-(*)}(\Omega))
\ \mbox{for any} \ 1 \leq q < \infty,
\eF
while
\bFormula{conv2}
\left\{ t \mapsto \intO{ n_\ep \vue \cdot \vc{H}^\perp [\varphi] } \right\}
\to 0 \ \mbox{in}\ L^q(0,T) \ \mbox{for any}\ 1 \leq q < \infty,\
\varphi \in \DC(\Omega).
\eF

\section{Compactness of the solenoidal part - proof of Theorem \ref{Tm1}}
\label{c}

In this section, we complete the proof of Theorem \ref{Tm1}. To begin,
we remark that relation (\ref{m1}) follows directly from (\ref{ub4}), while  (\ref{ub7a}) implies that
\bFormula{wc}
\vue \to \vc{U} \ \mbox{weakly in}\ L^2(0,T; W^{1,2}(\Omega;R^3),
\eF
at least for a subsequence as the case may be. Moreover, the vector field
$\vc{U}$ is solenoidal and satisfies the impermeability condition $\vc{U} \cdot \vc{n}|_{\partial \Omega} = 0$.

Next,
the uniform bounds (\ref{ub3}), (\ref{ub4}), together with the standard elliptic theory, yield
\bFormula{grad}
{\rm ess} \sup_{t \in (0,T)} \left\| \Grad  \left( \frac{\Phi_\ep(t, \cdot)}{\ep} \right) \right\|_{
W^{1,5/3}(K, R^3)} \leq c(K) \ \mbox{for any compact}\ K \subset \Ov{\Omega},
\eF
which, combined with (\ref{convv1}), yields
\bFormula{cc1}
\int_0^T \intO{ \frac{\nep - \Ov{n}}{\ep} \Grad \left( \frac{ \Phi_\ep}{\ep} \right)
\cdot \varphi }
\ \dt \to 0 \ \mbox{as}\ \ep \to 0 \ \mbox{for any} \ \varphi \in \DC([0,T)
\times \Ov{\Omega};R^3).
\eF

Taking $\varphi = \DC((0,T) \times \Omega;R^3)$, $\Div \varphi =0$, as a test function in the momentum equation(\ref{w2}) and making use of
(\ref{cc1}), we deduce that
\bFormula{c1}
\left\{ \vc{H}[\nep \vue] \right\}_{\ep > 0} \ \mbox{is precompact in}\
C_{\rm weak-(*)}([0,T]; L^2 + L^{5/4}(\Omega; R^3)).
\eF
Indeed the only quantity term reads
\[
\frac{1}{\ep^2} \nep \Grad \Phi_\ep = \frac{\Ov{n}}{\ep^2} \Grad \Phi_\ep +
\frac{\nep - \Ov{n}}{\ep} \Grad \frac{\Phi_\ep}{\ep},
\]
where the former term is a gradient, while the latter satisfies (\ref{cc1}).

Putting together (\ref{conv2}), (\ref{wc}), (\ref{c1}), with (\ref{m1}),
we deduce the desired conclusion
\bFormula{cccc}
\vue \to \vc{U} \ \mbox{(strongly) in} \ L^2((0,T) \times K; R^3)
\ \mbox{for any compact}\ K \subset \Omega.
\eF
With relations (\ref{cc1}), (\ref{cccc}) at hand, it is not difficult to perform the limit $\ep \to 0$ in the weak formulation of momentum equation
(\ref{si2}) to obtain (\ref{m7}).

We have proved Theorem \ref{Tm1}.

\section{Zero viscosity limit - proof of Theorem \ref{Tm2}}

\label{z}

Our ultimate goal is to prove Theorem \ref{Tm2}. The basic tool here is the
relative entropy inequality (\ref{w3}) satisfied by the  suitable weak solutions.
Taking $n_e = \nep$, $\vc{u} = \vue$, $r = \Ov{n}$ in
the rescaled variant of (\ref{w3}) we obtain
\bFormula{z1}
\intO{ \left( \frac{1}{2} n_\ep |\vue - \vc{U}|^2 +
\frac{1}{\ep^2}E(\nep, \Ov{n})
 + \frac{1}{2} \left|  \Grad  \left( \frac{\Phi_\ep }{\ep} \right)  \right|^2 \right) (s, \cdot)  }
\eF
\[
+ \int_0^s \intO{ \left[ \tn{S}_\ep (\Grad \vue) -
\tn{S}_\ep (\Grad \vc{U}) \right] : \Grad (\vue - \vc{U}) } \ \dt
+
\int_0^s \intO{ \frac{1}{\tau} \nep | \vue - \vc{U} |^2 } \ \dt
\]
\[
\leq \intO{ \left( \frac{1}{2} n_{0,\ep} |\vu_{0,\ep} - \vc{U}_{0, \ep} |^2 +
\frac{1}{\ep^2} E(n_{0,\ep}, \Ov{n}) + \frac{1}{2} \left| \Grad \left(
\frac{\Phi_{0,\ep}}{\ep} \right)  \right|^2  \right) } + \int_0^s {\cal R} \left(\nep,
\vue, \Ov{n}, \vc{U} \right) \ \dt,
\]
with
\bFormula{z2}
{\cal R}\left( n_e, \vue, \Ov{n}, \vc{U} \right) \equiv
\intO{  \nep  \Big( \partial_t \vc{U} + \vue \Grad \vc{U} \Big) \cdot
(\vc{U} - \vue )}
\eF
\[
+ \intO{\tn{S}_\ep(\Grad \vc{U}):\Grad (\vc{U}- \vue)  } +
\intO{ \frac{1}{\tau} \nep \vc{U} \cdot (\vc{U} - \vue ) }
- \frac{1}{\ep^2} \intO{ \nep \Grad \Phi_\ep \cdot \vc{U}}
\]
\[
- \frac{1}{\ep^2} \intO{ \Div \vc{U}  \Big( \nep \Big(
P(\nep) - P(\Ov{n}) \Big) - E (\nep, \Ov{n} ) \Big)  }.
\]

Furthermore, we take
\[
\vc{U} = \vc{U}_{\ep, \delta} = \vc{v} + \Grad \Psi_{\ep, \delta},
\]
where $\vc{v}$ is the (unique) solution of the damped Euler system
(\ref{euler1}-\ref{euler3}), emanating from the initial data
$\vc{v}_0 = \vc{H}[\vc{u}_0]$, and $\Grad \Psi_{\ep, \delta}$
mimicks the oscillatory part of the velocity field. Specifically, we take 
\bFormula{z3}
\ep \partial_t s_{\ep, \delta} + \Delta \Psi_{\ep, \delta} = 0,
\eF
\bFormula{z4}
\ep \partial_t \Grad \Psi_{\ep, \delta}  + p'(\Ov{n}) \Grad s_{\ep,
\delta} - \Ov{n} \Grad \Delta^{-1}_{N} s_{\ep, \delta} +
\frac{\ep}{\tau} \Grad \Psi_{\ep, \delta}= 0,\ \Grad \Psi_{\ep, \delta} \cdot \vc{n}|_{\partial \Omega} = 0,
\eF
which is nothing other than a slightly modified homogeneous part of
the acoustic system (\ref{Ab1}), (\ref{Ab2}). The initial data are
taken in the form
\bFormula{z4a}
s_{\ep, \delta}(0, \cdot) = \frac{1}{\Ov{n}} [N_{0,\ep}]_\delta, \ \Psi_{\ep, \delta}
(0, \cdot) = \left[ \Psi_{0,\ep} \right]_\delta,\ \mbox{with}\
\Grad \Psi_{0,\ep} = \vc{H}^\perp [\vu_{0,\ep}],
\eF
where the brackets $[ \cdot ]_{\delta}$ denote a suitable
regularization operator specified in Section \ref{decosc} below.

Keeping (\ref{z3} - \ref{z4}) in mind, we can rewrite the remainder (\ref{z2}) in the form
\bFormula{z6}
{\cal R}\left( n_e, \vue, \Ov{n}, \vc{U}_{\ep, \delta} \right)
\eF
\[
= \intO{  \nep  \Big( \partial_t \vc{v} + \vue \Grad \vc{v} + \frac{1}{\tau} \vc{v} \Big) \cdot
(\vc{U}_{\ep, \delta} - \vue )} + \intO{ \left( \nep \partial_t \Grad \Psi_{\ep, \delta}\cdot
\vc{v} + \nep \vue \cdot \Grad^2 \Psi_{\ep, \delta} \cdot ( \vc{U}_{\ep, \delta}-\vue) \right)}
\]
\[
+ \intO{\tn{S}_\ep(\Grad \vc{U}_{\ep, \delta}):\Grad (\vc{U}_{\ep, \delta} - \vue)  } -
\intO{ \Delta \Psi_{\ep, \delta} \left( \frac{P(\nep) - P(\Ov{n})}{\ep}
\frac{\nep - \Ov{n}}{\ep} - \frac{E(\nep, \Ov{n})}{\ep^2} \right) }
\]
\[
\intO{ \left( - \nep \vue \cdot \partial_t \Grad \Psi_{\ep, \delta} +
\frac{1}{2} \nep \partial_t | \Grad \Psi_{\ep, \delta} |^2 - \frac{\Ov{n}}{\ep}
\Delta \Psi_{\ep, \delta} \frac{ P(\nep) - P(\Ov{n}) }{\ep} \right) }
\]
\[
+
\intO{ \frac{1}{\tau}
 \nep \Grad \Psi_{\ep, \delta}  \cdot (\vc{U}_{\ep, \delta} - \vue ) }
- \frac{1}{\ep^2} \intO{ \nep \Grad \Phi_\ep \cdot \vc{U}_{\ep, \delta}  }.
\]

Moreover, we compute
\[
\intO{ \nep \vue \cdot \partial_t \Grad \Psi_{\ep, \delta} } =
p'(\Ov{n}) \intO{ N_\ep \partial_t s_{\ep, \delta} } - p'(\Ov{n})\left[
\intO{ N_\ep s_{\ep, \delta} } \right]_{t = 0}^{t = s}
\]
\[
- \Ov{n} \intO{ N_\ep \partial_t \Delta^{-1}_N [s_{\ep, \delta}] } +
\Ov{n} \left[ \intO{ N_\ep \Delta^{-1}_N [s_{\ep, \delta}] }
\right]_{t = 0}^{t = s} - \frac{1}{\tau} \intO{ \nep \vue \cdot
\Grad \Psi_{\ep, \delta}}
\]
and
\[
\intO{ \frac{1}{2} n_\ep \partial_t |\Grad \Psi_{\ep, \delta} |^2 } =
\intO{ \frac{1}{2} (n_\ep - \Ov{n})  \partial_t |\Grad \Psi_{\ep, \delta} |^2 } + \intO{ \frac{1}{2} \Ov{n} \partial_t |\Grad \Psi_{\ep, \delta} |^2 },
\]
where, by virtue of (\ref{z3}), (\ref{z4}),
\[
\intO{ \frac{1}{2} \Ov{n} \partial_t |\Grad \Psi_{\ep, \delta} |^2 }
= - \intO{ p'(n) \frac{\Ov{n}}{2} \partial_t | s_{\ep, \delta}|^2 }
+ \frac{\Ov{n}^2}{\ep} \intO{ \Grad \Psi_{\ep, \delta} \cdot \Grad
\Delta^{-1}_N [ s_{\ep, \delta} ] } - \frac{\Ov{n}}{\tau} \intO{
|\Grad \Psi_{\ep, \delta} |^2 }.
\]

Next, we have
\[
\frac{\Ov{n}}{\ep} \intO{ \Delta \Psi_{\ep,\delta} \frac{ P(\nep)
- P(\Ov{n}) } {\ep} }
\]
\[
= \frac{1}{\ep} \intO{ {p'(\Ov{n})} \Delta \Psi_{\ep,\delta} N_\ep }
+ \frac{\Ov{n} }{\ep^2} \intO{ \left( P(\nep) -
\frac{p'(\Ov{n})}{\Ov{n}} (\nep - \Ov{n} ) - P(\Ov{n}) \right)
\Delta \Psi_{\ep,\delta} },
\]
where
\[
\frac{1}{\ep} \intO{ \Delta \Psi_{\ep,\delta} N_\ep } = - \intO{
\partial_t s_{\ep,\delta} N_\ep }.
\]

Finally,
\[
\frac{1}{\ep^2} \intO{ \nep \Grad \Phi_\ep \cdot \vc{U}_{\ep,\delta}
} = \intO{ N_\ep \Grad \left( \frac{\Phi_\ep} {\ep} \right) \cdot
\vc{U}_{\ep, \delta} } + \frac{\Ov{n}}{\ep} \intO{\Grad \left(
\frac{\Phi_\ep} {\ep} \right) \cdot \Grad \Psi_{\ep, \delta} }
\]
\[
=\intO{ N_\ep \Grad \left( \frac{\Phi_\ep} {\ep} \right) \cdot
\vc{U}_{\ep, \delta} } + \Ov{n} \intO{ N_\ep \partial_t
\Delta^{-1}_N [ s_{\ep, \delta}] }.
\]

Summing up the previous considerations we may infer that
\bFormula{z9}
\intO{ \left( \frac{1}{2} n_\ep |\vue - \vc{v} - \Grad \Psi_{\ep,
\delta} |^2 + \frac{1}{\ep^2}E(\nep, \Ov{n})
 + \frac{1}{2} \left|  \Grad  \left( \frac{\Phi_\ep }{\ep} \right)  \right|^2 \right) (s, \cdot)  }
\eF
\[
+ \intO{\left( p'(\Ov{n}) \frac{\Ov{n}}{2} |s_{\ep,\delta}|^2 -
p'(\Ov{n}) N_\ep s_{\ep,\delta} + \Ov{n} N_\ep \Delta^{-1}_N [
s_{\ep,\delta}]  + \frac{\Ov{n}^2}{2} \left| (-\Delta_n)^{-1/2}
[s_{\ep, \delta}] \right|^2  \right) (s, \cdot) }
\]
\[
+ \frac{\mu_\ep}{2} \int_0^s \intO{ \left| \Grad (\vue -
\vc{U}_{\ep,\delta} ) + \Grad^t (\vue - \vc{U}_{\ep, \delta}) -
\frac{2}{3} \Div (\vue - \vc{U}_{\ep, \delta}) \tn{I} \right|^2 } \
\dt + \int_0^s \intO{ \frac{1}{\tau} \nep | \vue - \vc{U}_{\ep,
\delta} |^2 } \ \dt
\]
\[
\leq \intO{ \left(  + \frac{1}{2} n_{0,\ep} \Big| \vu_{0,\ep} -
\vc{H}[\vu_0] - \Grad [\Psi_{0, \ep}]_\delta \Big|^2 +
\frac{1}{\ep^2} E(n_{0,\ep}, \Ov{n}) + \frac{1}{2} \left| \Grad
\left( \frac{\Phi_{0,\ep}}{\ep} \right) \right|^2  \right) } +
\]
\[
+ \intO{\left( p'(\Ov{n}) \frac{\Ov{n}}{2} |[N_{0, \ep}]_{\delta}|^2
- p'(\Ov{n}) N_{0,\ep} [N_{0, \ep}]_{\delta} + \Ov{n} N_{0,\ep}
\Delta^{-1}_N [ N_{0,\ep}]_\delta + \frac{\Ov{n}^2}{2} \left|
(-\Delta_n)^{-1/2} [N_{0,\ep}]_\delta \right|^2 \right) }
\]
\[
+ \int_0^s Q_{\ep, \delta} \ \dt,
\]
where
\[
Q_{\ep, \delta} =
\]
\[
= \intO{  \nep  \Big( \partial_t \vc{v} + \vue \Grad \vc{v} +
\frac{1}{\tau} \vc{v} \Big) \cdot (\vc{U}_{\ep, \delta} - \vue )} +
\intO{ \left( \nep \partial_t \Grad \Psi_{\ep, \delta}\cdot \vc{v} +
\nep \vue \cdot \Grad^2 \Psi_{\ep, \delta} \cdot (\vc{U}_{\ep, \delta} - \vue ) \right)}
\]
\[
+ \frac{1}{\tau} \intO{ (\nep - \Ov{n})  \vc{U}_{\ep, \delta} \cdot
\Grad \Psi_{\ep, \delta} } +\frac{1}{\tau}\intO{\nep(\vue-\vc{U}_{\ep,\delta})\cdot\Grad\Psi_{\ep,\delta}} + \frac{1}{2} \intO{ (\nep - \Ov{n})
\partial_t | \Grad \Psi_{\ep, \delta}|^2 }
\]
\[
+ \intO{\tn{S}_\ep(\Grad \vc{U}_{\ep, \delta}):\Grad (\vc{U}_{\ep,
\delta} - \vue)  } - \intO{ \Delta \Psi_{\ep, \delta} \left(
\frac{P(\nep) - P(\Ov{n})}{\ep} \frac{\nep - \Ov{n}}{\ep} -
\frac{E(\nep, \Ov{n})}{\ep^2} \right) }
\]
\[
- \frac{\Ov{n}}{\ep^2} \intO{ \left( P(\nep) -
\frac{p'(\Ov{n})}{\Ov{n}} (\nep - \Ov{n}) - P (\Ov{n}) \right)
\Delta \Psi_{\ep, \delta} } - \intO{ N_\ep \Grad \left( \frac{
\Phi_\ep}{\ep} \right) \cdot \vc{U}_{\ep, \delta} }.
\]

\subsection{Dispersive estimates of the oscillatory component}

\label{decosc}

Our goal is to show that solutions $s_{\ep,\delta}$, $\Psi_{\ep, \delta}$ of the 
homogeneous ``acoustic'' equation (\ref{z3}), (\ref{z4}) decay to zero in the $L^\infty$ norm as $\ep \to 0$ for any positive time $t$. 
To this end, we start with the total energy balance 
\bFormula{ener}
\frac{{\rm d}}{{\rm d}t} \intO{ \left( |\Grad \Psi_{\ep, \delta}|^2 + p'(\Ov{n}) |s_{\ep, \delta}|^2 + 
\Ov{n} \left| (-\Delta_N)^{-1/2}[s_{\ep, \delta}] \right|^2 \right) }
+ \frac{2}{\tau} \intO{ |\Grad \Psi_{\ep, \delta} |^2 } = 0
\eF
yielding, in particular, existence and uniqueness of (weak) solutions to 
problem (\ref{z3}), (\ref{z4}) provided the initial data are smooth and decay sufficiently fast for $|x| \to \infty$. 

Taking advantage of the special geometry of waveguides, we consider 
the functions $w_k(z)$, $z \in B$ - the eigenfunctions of the Neumann 
Laplacian $-\Delta_{N,B}$ in the (bounded) domain $B \subset R^{3 - L}$: 
\[
- \Delta_{N,B} w_k = \lambda_k w_k \ \mbox{in} \ B, 
\ \nabla_z w_k \cdot \vc{n}_z |_{\partial B} = 0,\ \lambda_0 = 0 < 
\lambda_1 \leq \dots \leq \lambda_k \leq \dots , \ k = 0,1, \dots.
\]

The smoothing operators $[g]_\delta$, $g = g(x),\ x =[y,z]$ are defined as 
\bFormula{smooth}
[g]_\delta (y,z) = \sum_{0 \leq k < 1/\delta} \kappa_\delta (y) * \Big( \psi_\delta (y) A_k [g](y)\Big) w_k(z),
\eF
where 
\[
A_k[g] (y) = \frac{1}{|B|} \int_B g(y,z) w_k(z) \ {\rm d}z,
\]
$\psi_\delta \in \DC(R^L)$ is a cut-off function, 
\[
0 \leq \psi_\delta \leq 1, \ \psi_\delta (y) = \left\{ 
\begin{array}{l} 1 \ \mbox{for}\ |y| < 1/\delta , \\ \\ 
0 \ \mbox{for}\ |y| > 2/\delta, \end{array} \right.
\]
and $\kappa_\delta$ is a family of standard regularizing kernels in the 
$y-$variable.

A short inspection of (\ref{z3}), (\ref{z4}) yields
\bFormula{formulk}
s_{\ep, \delta} (t,x) = \exp \left( \frac{t}{2\tau} \right) \tilde s_{\ep, \delta} \left( \frac{t}{\ep}, x \right),
\eF
where $\tilde s_{\ep, \delta}$ is the unique solution of the Klein-Gordon 
equation 
\[
\partial^2_{t,t} \tilde s_{\ep, \delta} - p'(\Ov{n}) \Delta \tilde s_{\ep, \delta} + \left( \Ov{n} - \frac{1}{4} \frac{\ep^2}{\tau_2} \right) \tilde s_{\ep, \delta} = 0,\ \Grad \tilde s_{\ep, \delta} \cdot \vc{n}|_{\partial \Omega} = 0, 
\]
emanating from the initial data 
\bFormula{iniD}
\tilde s_{\ep, \delta} (0, \cdot)  = \frac{1}{\Ov{n}} [N_{0,\ep}]_{\delta},\ 
\partial_t \tilde s_{\ep, \delta} (0, \cdot)  = - \Delta [\Psi_{0,\ep}]_\delta.
\eF

Consequently, thanks to the specific choice of the smoothing operators 
(\ref{smooth}), solutions $\tilde s_{\ep, \delta}$ take the form 
\[
\tilde s_{\ep, \delta}(t,x) = \sum_{0 \leq k \leq 1/\delta} 
S_{k, \ep, \delta} (t,y) w_k(z), 
\]
where $S_k(t, \cdot)$ solve the Klein-Gordon equation
\bFormula{KG}
\partial^2_{t,t} S_{k,\ep, \delta} - p'(\Ov{n}) \Delta_y S_{k,\ep, \delta} + \left( \Ov{n} - \frac{1}{4} \frac{\ep^2}{\tau^2} + p'(\Ov{n}) \lambda_k \right)  S_{k,\ep, \delta} = 0
\eF
for $y$ belonging to the ``flat'' space $R^L$, and with the initial data uniquely determined through (\ref{iniD}). Thus, employing the standard 
$L^1-L^\infty$ estimates for the Klein-Gordon equation \ref{KG} (see for instance Lesky and Racke \cite[Lemma 2.4]{LesRa2}), we have
\[
\left\| S_{k,\ep, \delta}(t, \cdot) \right\|_{L^\infty(R^L)} \leq 
\frac{c(k)}{ (1 + t)^{L/2}} \Big( \left\| S_{k,\ep, \delta}(0,\cdot) \right\|
_{W^{K,1}(R^L)} + \left\| \partial_t S_{k,\ep, \delta}(0,\cdot) \right\|
_{W^{K-1,1}(R^L)} \Big) , \ K = \left[ \frac{L + 3}{2} \right]. 
\]
Going back to (\ref{formulk}) we may infer that
\bFormula{DECAY1}
\| s_{\ep, \delta}(t, \cdot) \|_{L^\infty(\Omega)} \leq \omega (t_0, \ep, \delta),\ t \in [t_0, T], \ t_0 > 0,   
\eF
and, using (\ref{z3}),
\bFormula{DECAY2}
\| \Delta \Psi_{\ep, \delta}(t, \cdot) \|_{L^\infty(\Omega)} \leq \omega (t_0, \ep, \delta),\ t \in [t_0, T], \ t_0 > 0,   
\eF
where $\omega(t_0, \ep, \delta) \to 0$ if $\ep \to 0$ for any fixed $t_0> 0$, $\delta > 0$.

Finally, we claim the standard energy bounds
\bFormula{dec1}
\| \Grad \Psi_{\ep, \delta} \|_{L^\infty(0,T;W^{m,2}(\Omega;R^3))} \leq c (m, \delta), 
\eF
\bFormula{dec2}
\| s_{\ep, \delta} \|_{L^\infty(0,T;W^{m,2}(\Omega;R^3))} \leq c(m,  \delta), \ m = 0,1,\dots, 
\eF 
where the constants are independent of $\ep$ for any fixed $\delta > 0$.

\subsection{Asymptotic limit $\ep \to 0$}

Our next goal is to let $\ep \to 0$ in (\ref{z9}), and, in particular, in the remainder $Q_{\ep, \delta}$. 

\begin{enumerate}
\item 

We have
\[
\intO{  \nep  \Big( \partial_t \vc{v} + \vue \Grad \vc{v} +
\frac{1}{\tau} \vc{v} \Big) \cdot (\vc{U}_{\ep, \delta} - \vue )}
\]
\[
=
\intO{  \nep  \Big( \partial_t \vc{v} + \vc{U}_{\ep, \delta} \Grad \vc{v} +
\frac{1}{\tau} \vc{v} \Big) \cdot (\vc{U}_{\ep, \delta} - \vue )} 
+
\intO{  \nep  (\vue - \vc{U}_{\ep, \delta}) \cdot \Grad \vc{v} \cdot (\vc{U}_{\ep, \delta} - \vue )},
\]
where 
\[
\intO{  \nep  \Big( \partial_t \vc{v} + \vc{U}_{\ep, \delta} \Grad \vc{v} +
\frac{1}{\tau} \vc{v} \Big) \cdot (\vc{U}_{\ep, \delta} - \vue )} 
\]
\[
=\intO{  \nep  \Grad \Pi \cdot (\vue - \vc{U}_{\ep, \delta})} + 
\intO{ \nep \Grad \Psi_{\ep, \delta} \cdot \Grad \vc{v} \cdot (\vc{U}_{\ep, \delta} - \vue ) },
\]
and
\[
\intO{  \nep  \Grad \Pi \cdot (\vue - \vc{U}_{\ep, \delta})} = \intO{  \nep  \Grad \Pi \cdot \vue }
- \ep \intO{ N_\ep \Grad \Pi \cdot (\vc{v} + \Grad \Psi_{\ep, \delta}) }
+ \Ov{n} \intO{ \Grad \Pi \cdot \Grad \Psi_{\ep, \delta} }.
\]

Since 
\[
\nep \vue = [ \sqrt{ \nep } ]_{\rm ess} \sqrt{\nep} \vue + 
[\sqrt{\nep}]_{\rm res} \sqrt{\nep} \vue, 
\]
where, by virtue of estimates (\ref{ub2}), (\ref{ub4}), 
\[
[\sqrt{\nep}]_{\rm ess} \sqrt{\nep} \vue \to \Ov{n} \vc{U} \ \mbox{weakly-(*) in}\ 
L^\infty (0,T; L^2(\Omega;R^3)),\ \Div \vc{U} = 0, 
\]
while
\[
[\sqrt{\nep}]_{\rm res} \sqrt{\nep} \vue \to 0 \ \mbox{in}\ L^\infty(0,T; L^{5/4}(\Omega)),
\]
we get 
\[
{\rm ess} \sup_{t \in (0,T_{\rm loc})} \left| \intO{  \nep  \Grad \Pi \cdot \vue } \right| \to 0 \ \mbox{for}\ \ep \to 0.
\]

Similarly, we use (\ref{dec1}) to observe that 
\[
\left\{
t \mapsto \intO{ \Grad \Pi \cdot \Grad \Psi_{\ep, \delta} }
\right\} \to 0 \ \mbox{in}\ L^2(0,T_{\rm loc}) \ \mbox{as} \ \ep \to 0
\]
for any fixed $\delta > 0$. 

Thus we conclude that 
\bFormula{es1}
\left|
\intO{  \nep  \Big( \partial_t \vc{v} + \vue \Grad \vc{v} +
\frac{1}{\tau} \vc{v} \Big) \cdot (\vc{U}_{\ep, \delta} - \vue )}
\right| \leq c \intO{ n_\ep | \vue - \vc{U}_{\ep, \delta} |^2 } + 
h^1_{\ep, \delta}, 
\eF
where 
\bFormula{es2}
h^1_{\ep, \delta} \to 0 \ \mbox{in} \ L^2(0,T_{\rm loc}) \ \mbox{as}\ \ep \to 0.
\eF

\item 

Taking advantage of the fact that $\Div \vc{v} = 0$ we can write 
\[
\intO{ \nep \partial_t \Grad \Psi_{\ep, \delta} \cdot \vc{v} } =
\ep \intO{ N_\ep \partial_t \Grad \Psi_{\ep, \delta} \cdot \vc{v} }, 
\]
where $\ep \partial_t \Grad \Psi_{\ep, \delta}$ can be expressed by means of equation (\ref{z4}). Using (\ref{dec1}), (\ref{dec2}) we conclude that 
\bFormula{s3}
\left|
\intO{ \nep \partial_t \Grad \Psi_{\ep, \delta} \cdot \vc{v} } 
\right| = h^2_{\ep, \delta}, 
\eF
with
\bFormula{s4}
h^2_{\ep, \delta} \to 0 \ \mbox{in}\ L^2(0,T_{\rm loc}) 
\ \mbox{as}\ \ep \to 0.
\eF

\item Using (\ref{DECAY2}), (\ref{dec1}), we 
show that
\bFormula{s5}
\left| \intO{\nep \vue \cdot \Grad^2 \Psi_{\ep, \delta} \cdot (\vue -
\vc{U}_{\ep, \delta})} \right| +\left|\intO{\nep(\vue-\vc{U}_{\ep,\delta})\cdot\Grad\Psi_{\ep,\delta}}\right|
\eF
\[
+ \left| \intO{ (\nep - \Ov{n})  \vc{U}_{\ep, \delta} \cdot
\Grad \Psi_{\ep, \delta} } \right| + \left| \intO{ (\nep - \Ov{n})
\partial_t | \Grad \Psi_{\ep, \delta}|^2 } \right| \leq h^3_{\ep, \delta},
\]
with 
\bFormula{s6}
h^3_{\ep, \delta} \to 0 \ \mbox{in}\ L^2(0,T_{\rm loc}) \ \mbox{as}\ \ep \to 0.
\eF

\item

Now, 
\bFormula{s7}
\intO{ \tn{S}_\ep (\Grad \vc{U}_{\ep, \delta}) : \Grad (\vc{U}_{\ep, \delta} - \vue ) } 
\eF
\[
\leq \frac{\mu_\ep}{2} \intO{ \left| \Grad (\vue -
\vc{U}_{\ep,\delta} ) + \Grad^t (\vue - \vc{U}_{\ep, \delta}) -
\frac{2}{3} \Div (\vue - \vc{U}_{\ep, \delta}) \tn{I} \right|^2 } 
\]
\[
+ c \mu_\ep \intO{ |\Grad \vc{U}_{\ep, \delta} |^2 }.
\]

\item 

Next, in accordance with (\ref{ub4}) and (\ref{DECAY2}), (\ref{dec1}), 
\bFormula{s8}
\left| \intO{ \Delta \Psi_{\ep, \delta} \left(
\frac{P(\nep) - P(\Ov{n})}{\ep} \frac{\nep - \Ov{n}}{\ep} -
\frac{E(\nep, \Ov{n})}{\ep^2} \right) } \right| 
\eF
\[
+
\left|\frac{\Ov{n}}{\ep^2} \intO{ \left( P(\nep) -
\frac{p'(\Ov{n})}{\Ov{n}} (\nep - \Ov{n}) - P (\Ov{n}) \right)
\Delta \Psi_{\ep, \delta} }\right| \leq h^4_{\ep, \delta} \to 0 \ \mbox{in}\ 
L^2(0,T_{\rm loc}) \ \mbox{as}\ \ep \to 0.
\]

\item 

Finally, using (\ref{convv1}), (\ref{grad}), and (\ref{cc1}), we infer that
\bFormula{s9}
\left|
\intO{ N_\ep \Grad \left( \frac{
\Phi_\ep}{\ep} \right) \cdot \vc{U}_{\ep, \delta} } \right| 
= h^5_{\ep, \delta} \to 0 \ \mbox{in}\ L^2(0,T) \ \mbox{as}\ \ep \to 0.
\eF

\end{enumerate}

Using estimates (\ref{es1} - \ref{s9}) in (\ref{z9}) we conclude that
\bFormula{gron}
\intO{ \left( \frac{1}{2} n_\ep |\vue - \vc{v} - \Grad \Psi_{\ep,
\delta} |^2 \right) (s, \cdot) }
\eF
\[
+ \intO{ \left(  \frac{1}{\ep^2}E(\nep, \Ov{n})   -
p'(\Ov{n}) N_\ep s_{\ep,\delta} + p'(\Ov{n}) \frac{\Ov{n}}{2} |s_{\ep,\delta}|^2 \right) (s, \cdot) }
\]
\[
+ \intO{ \left( \frac{1}{2} \left|  \Grad  \left( \frac{\Phi_\ep }{\ep} \right) \right|^2
+ \Ov{n} N_\ep \Delta^{-1}_N [
s_{\ep,\delta}]  + \frac{\Ov{n}^2}{2} \left| (-\Delta_n)^{-1/2}
[s_{\ep, \delta}]   \right|^2 \right) (s, \cdot) }
\]
\[
\leq \intO{ \left( \frac{1}{2} n_{0,\ep} |\vu_{0,\ep} - \vc{H}[\vu_0] - \Grad [\Psi_{0,\ep}]_\delta
 |^2 \right) }
\]
\[
+ \intO{ \left(  \frac{1}{\ep^2}E(n_{0,\ep}, \Ov{n})   -
p'(\Ov{n}) N_{0,\ep} [N_{0,\ep}]_\delta + p'(\Ov{n}) \frac{\Ov{n}}{2} \Big| [N_{0,\ep}]_\delta \Big|^2 \right) }
\]
\[
+ \intO{ \left( \frac{1}{2} \left|  \Grad  \left( \frac{\Phi_{0,\ep}}{\ep} \right) \right|^2
+ \Ov{n} N_{0,\ep} \Delta^{-1}_N 
[N_{0,\ep}]_\delta  + \frac{\Ov{n}^2}{2} \left| (-\Delta_n)^{-1/2}
[N_{0,\ep}]_\delta  \right|^2 \right)  }
\]
\[
+c \int_0^s \intO{ \left( \frac{1}{2} n_\ep |\vue - \vc{v} - \Grad \Psi_{\ep,
\delta} |^2 \right)  } \ {\rm d}t + \int_0^s h^6_{\ep, \delta} \ \dt
\]
where 
\bFormula{gron1}
h^6_{\ep, \delta} \to 0 \ \mbox{in} \ L^2(0,T_{\rm loc}) \ \mbox{as} \ \ep \to 0.
\eF

Now, we claim that 
\[
{\rm ess} \sup_{t \in (0,T)} \left\| \sqrt{ \frac{ E(\nep, \Ov{n})}{\ep^2} } - \sqrt{ \frac{p'(\Ov{n})}{ 2 \Ov{n}}}
N_\ep \right\|_{L^{4/3}(\Omega)} \to 0 \ \mbox{as}\ \ep \to 0.
\]
Indeed we have
\[
 \sqrt{ \frac{ E(\nep, \Ov{n}) }{\ep^2} } - \sqrt{ \frac{p'(\Ov{n})}{2 \Ov{n}}} N_\ep
\]
\[
= \left[ \left( \frac{ H(\nep) - H'(\Ov{n}) (\nep - \Ov{n}) -
H(\Ov{n}) }{\ep^2} \right)^{1/2} - \left( \frac{1}{2}
H''(\Ov{n}) \frac{(\nep - \Ov{n})^2}{\ep^2} \right)^{1/2}
\right]_{\rm ess}
\]
\[
+\left[ \left( \frac{ H(\nep) - H'(\Ov{n}) (\nep - \Ov{n}) -
H(\Ov{n}) }{\ep^2} \right)^{1/2} - \left( \frac{1}{2}
H''(\Ov{n}) \frac{(\nep - \Ov{n})^2}{\ep^2} \right)^{1/2}
\right]_{\rm res},
\]
where
\[
\left[ \left( \frac{ H(\nep) - H'(\Ov{n}) (\nep - \Ov{n}) -
H(\Ov{n}) }{\ep^2} \right)^{1/2} - \left( \frac{1}{2}
H''(\Ov{n}) \frac{(\nep - \Ov{n})^2}{\ep^2} \right)^{1/2}
\right]_{\rm res} \to 0
\]
\[
\mbox{in}\ L^\infty(0,T; L^q(\Omega)), \ 1 \leq q < 5/3,
\]
while
\[
\left| \left[ \left( \frac{ H(\nep) - H'(\Ov{n}) (\nep - \Ov{n})
- H(\Ov{n}) }{\ep^2} \right)^{1/2} - \left( \frac{1}{2}
H''(\Ov{n}) \frac{(\nep - \Ov{n})^2}{\ep^2} \right)^{1/2}
\right]_{\rm ess} \right|
\]
\[
\leq \sqrt{\ep} \left[ \sqrt{ H'''(\xi) \frac{ |\nep - \Ov{n}|^3 }
{\ep^3} } \right]_{\rm ess} \ \mbox{for a certain}\ \xi \in
[\Ov{n}/2, 2 \Ov{n}].
\]

Consequently, relation (\ref{gron}) can be written in the form 
\bFormula{gron2}
\intO{ \left( \frac{1}{2} n_\ep |\vue - \vc{v} - \Grad \Psi_{\ep,
\delta} |^2 \right) (s, \cdot) }
\eF
\[
+  \intO{ \left(  \sqrt{ \frac{E(\nep, \Ov{n})} {\ep^2
 }} -
\sqrt{ \frac{ p'(\Ov{n}) \Ov{n} }{2} } s_{\ep,\delta} \right)^2 (s, \cdot) }
+ \intO{ \frac{1}{2} \left| (-\Delta)^{-1/2} \left[ s_{\ep, \delta} - N_\ep \right]  \right|^2  (s, \cdot) }
\]
\[
\leq \intO{ \left( \frac{1}{2} n_{0,\ep} |\vu_{0,\ep} - \vc{H}[\vu_0] - \Grad [\Psi_{0,\ep}]_\delta
 |^2 \right) }
\]
\[
+ \intO{ \left(  \frac{1}{\ep^2}E(n_{0,\ep}, \Ov{n})   -
p'(\Ov{n}) N_{0,\ep} [N_{0,\ep}]_\delta + p'(\Ov{n}) \frac{\Ov{n}}{2} \Big| [N_{0,\ep}]_\delta \Big|^2 \right) }
\]
\[
+ \intO{ \left( \frac{1}{2} \left|  \Grad  \left( \frac{\Phi_{0,\ep}}{\ep} \right) \right|^2
+ \Ov{n} N_{0,\ep} \Delta^{-1}_N 
[N_{0,\ep}]_\delta  + \frac{\Ov{n}^2}{2} \left| (-\Delta_n)^{-1/2}
[N_{0,\ep}]_\delta  \right|^2 \right)  }
\]
\[
+c \int_0^s \intO{ \left( \frac{1}{2} n_\ep |\vue - \vc{v} - \Grad \Psi_{\ep,
\delta} |^2 \right)  } \ {\rm d}t + \int_0^s h^7_{\ep, \delta} \ \dt
\]
with 
\[
h^7_{\ep, \delta} \to 0 \ \mbox{in}\ L^2(0,T_{\rm loc})\ \mbox{as}\ \ep \to 0.
\]

Applying Gronwall's lemma we therefore get
\[
\intO{ \left( \frac{1}{2} n_\ep |\vue - \vc{v} - \Grad \Psi_{\ep,
\delta} |^2 \right) (s, \cdot) }
\]
\[
+ \intO{ \left(  \sqrt{ \frac{E(\nep, \Ov{n})} {\ep^2
 }} -
\sqrt{ \frac{ p'(\Ov{n}) \Ov{n} }{2} } 
s_{\ep,\delta} \right)^2 (s, \cdot) }
+ \intO{ \frac{1}{2} \left| (-\Delta)^{-1/2} \left[ s_{\ep, \delta} - N_\ep \right]  \right|^2  (s, \cdot) }
\]
\[
\leq R_{\ep,\delta}(s)+ \left( C\int_0^s R_{\ep,\delta}(t)
e^{- Ct}{\rm d}t\right) e^{Cs},
\]
where 
\[
R_{\ep, \delta}(s) =  \intO{ \left( \frac{1}{2} n_{0,\ep} |\vu_{0,\ep} - \vc{H}[\vu_0] - \Grad [\Psi_{0,\ep}]_\delta
 |^2 \right) }
\]
\[
+ \intO{ \left(  \frac{1}{\ep^2}E(n_{0,\ep}, \Ov{n})   -
p'(\Ov{n}) N_{0,\ep} [N_{0,\ep}]_\delta + p'(\Ov{n}) \frac{\Ov{n}}{2} \Big| [N_{0,\ep}]_\delta \Big|^2 \right) }
\]
\[
+ \intO{ \left( \frac{1}{2} \left|  \Grad  \left( \frac{\Phi_{0,\ep}}{\ep} \right) \right|^2
+ \Ov{n} N_{0,\ep} \Delta^{-1}_N 
[N_{0,\ep}]_\delta  + \frac{\Ov{n}^2}{2} \left| (-\Delta_n)^{-1/2}
[N_{0,\ep}]_\delta  \right|^2 \right)  } + \int_0^s h^7_{\ep, \delta} \ \dt.
\]
Thus, letting $\ep \to 0$ we obtain 
\bFormula{konec1}
{\rm ess} \sup_{\ep \to 0} \sup_{s \in (0, T_{\rm loc})}\intO{ \left( \frac{1}{2} n_\ep |\vue - \vc{v} - \Grad \Psi_{\ep,
\delta} |^2 \right) (s, \cdot) } \le \chi(\delta),
\eF
where the function $\chi$ is determined in terms of the initial data, and $\chi(\delta) \to 0$ as $\delta \to 0$.

\subsection{Asymptotic limit $\delta \to 0$}

Letting $\delta \to 0$ in (\ref{konec1}) we may infer that 
\bFormula{konec2}
\vc{ H}[\sqrt\nep\,\vue]\to\sqrt{\overline n}\,\vc v\;\mbox{in
$L^2(0,T_{\rm loc};L^2(K;R^3))$},
\eF
\bFormula{konec3}
\vc{H}^\perp[\sqrt\nep\,\vue]\to 0 \;\mbox{in $L^2(0,T_{\rm loc};L^2(K;R^3))$}
\eF
for any compact $K \subset \Omega$.

Relations (\ref{konec2}), (\ref{konec3}) complete the proof of Theorem \ref{Tm2}.

\def\ocirc#1{\ifmmode\setbox0=\hbox{$#1$}\dimen0=\ht0 \advance\dimen0
  by1pt\rlap{\hbox to\wd0{\hss\raise\dimen0
  \hbox{\hskip.2em$\scriptscriptstyle\circ$}\hss}}#1\else {\accent"17 #1}\fi}

%\bibliography{citace}
%\bibliographystyle{plain}

\end{document}